\documentclass[12pt,reqno]{article}
\usepackage{amssymb, amsthm, amsmath, amsfonts, stmaryrd, latexsym, wasysym}
\usepackage{enumitem}
\usepackage[v2]{xy}
\xyoption{curve}

\usepackage{amssymb}
\usepackage{graphicx}
\usepackage{amscd}
\usepackage{color}

\usepackage{fullpage}
\usepackage{setspace}

\usepackage{fullpage}
\usepackage{float}
\usepackage{graphicx,amsmath,amssymb}
\usepackage{amsthm}
\usepackage{amsfonts}
\usepackage{latexsym}
\usepackage{epsf}

\definecolor{grau}{rgb}{0.3,0.3,0.2}

\definecolor{rot}{rgb}{0.3,0.3,0.2}
\definecolor{gruen}{rgb}{0.3,0.3,0.2}
\definecolor{blau}{rgb}{0.3,0.3,0.2}

\def\st 				{\ |\ }
\def\na				{\vskip 0.12cm}
\def\gra				{\vskip 0.36cm}
\def\ne 				{\noindent}

\def\tabulaturna     {\vskip 0.12cm}
\def\tabulaturkla    {\vskip 0.06cm}

\def\qed				{\hfill $\Box$}
\def\qee				{\hfill $\Box$}
\def\ns				{}

\def\tabulaturna     {\vskip 0.12cm}
\def\tabulaturkla    {\vskip 0.12cm}
\def\loss            {\begin{tabular}{p{0.30in}p{5.95in}}}

\begin{document}

\color{grau}{

\begin{center}

\vskip 1cm{\large\bf Solvable Hypergroups and a Generalization of Hall's}

\vskip .32cm{\large\bf Theorems on Finite Solvable Groups to Association Schemes}

\vskip .8cm
\tiny
Andrey Vasil'ev \\
Sobolev Institute of Mathematics \\
4 Acad.\ Koptyug avenue \\
Novosibirsk, 630090, Russia \\
and \\
Novosibirsk State University \\
1 Pirogova street \\
Novosibirsk, 630090, Russia \\

\vskip .45cm

Paul-Hermann Zieschang \\
School of Mathematical and Statistical Sciences \\
University of Texas Rio Grande Valley \\
Edinburg, TX 78539, U.\ S.\ A. \\
\end{center}

\vskip 0.5cm

\begin{abstract}

\ne We generalize Philip Hall's celebrated theorems on finite solvable groups to scheme theory. Our result is based on a series of results on hypergroups.

\end{abstract}

\vskip .5cm

\ns\centerline{\bf 1. Introduction} \gra\na

\ne The concept of an association scheme provides a far-reaching and meaningful generalization of the concept of a group. A number of important results on finite groups have already found generalizations to association schemes. In {\color{gruen} [6]}, for instance, analogues of the homomorphism theorem, the isomorphism theorems, and the Jordan-H\"{o}lder Theorem for finite groups were proved for association schemes. In {\color{gruen} [4]}, a scheme theoretic generalization of Sylow's Theorems was proved, and in {\color{gruen} [1]}, a corresponding generalization of the Schur-Zassenhaus Theorem on finite groups was given. \na

\ne It is the purpose of the present note to present a generalization of Philip Hall's celebrated theorem {\color{gruen} [2; Theorem]} on finite solvable groups to scheme theory. We first introduce notation and terminology needed for the statement of our result. We begin with a review of the definition of an association scheme. \na

\ns\ne Let $X$ be a finite set. We write $1_{X}$ to denote the set of all pairs $(x,x)$ with $x\in X$. For each subset $r$ of the cartesian product $X\times X$, we define $r^{*}$ to be the set of all pairs $(y,z)$ with $(z,y)\in r$. Whenever $x$ is an element in $X$ and $r$ a subset of $X\times X$, we write $xr$ to denote the set of all elements $y$ in $X$ with $(x,y)\in r$. \na

\ne Let $S$ be a partition of $X\times X$ with $1_{X}\in S$, and assume that $s^{*}\in S$ for each element $s$ in $S$. The set $S$ is called an {\color{blau} \em association scheme} or simply a {\color{blau} \em scheme on} $X$ if, for any three elements $p$, $q$, and $r$ in $S$, there exists an integer $a_{pqr}$ such that $|yp\cap zq^{*}|=a_{pqr}$ for any two elements $y$ in $X$ and $z$ in $yr$.\footnote{\color{grau} We note that the definition of a scheme which we use in the present note differs from the more general definition given in {\color{gruen} [8]}. In fact, the schemes which we consider in this note are exactly the schemes in the sense of {\color{gruen} [8]} which are defined on finite sets.} \na

\ns\ne Let $S$ be an association scheme, and let $P$ and $Q$ be subsets of $S$. We define $PQ$ to be the set of all elements $s$ in $S$ such that there exist elements $p$ in $P$ and $q$ in $Q$ with $a_{pqs}\neq 0$. The set $PQ$ is called the {\color{blau} \em complex product} of $P$ and $Q$, and the map which assigns to any two subsets of $S$ its complex product will be called the {\color{blau} \em complex multiplication in} $S$. If $P$ contains an element $p$ with $P=\{p\}$, we write $pQ$ instead of $PQ$. Similarly, if $Q$ contains an element $q$ with $Q=\{q\}$, we write $Pq$ instead of $PQ$. Finally, if $P$ contains an element $p$ with $P=\{p\}$ and $Q$ contains an element $q$ with $Q=\{q\}$, we write $pq$ instead of $PQ$. \na

\ne A non-empty subset $R$ of an association scheme is called {\color{blau} \em closed} if $p^{*}q\subseteq R$ for any two elements $p$ and $q$ in $R$. \na

\ne Let $S$ be a scheme on a set $X$. Instead of $1_{X}$ we write $1$. For each element $s$ in $S$, we abbreviate $n_{s}:=a_{ss^{*}1}$, and we call $n_{s}$ the {\color{blau} \em valency} of $s$. For each closed subset $T$ of $S$, the sum of the positive integers $n_{t}$ with $t\in T$ will be denoted by $n_{T}$, and we call $n_{T}$ the {\color{blau} \em valency} of $T$. \na

\ne A closed subset of an association scheme is called {\color{blau} \em thin} if all of its elements have valency $1$. \na

\ns\ne There is an obvious way to identify groups with thin association schemes. This correspondence was established in {\color{gruen} [8; Section 5.5]} and is called the {\color{blau} \em group correspondence}. \na

\ne Let $\pi$ be a set of prime numbers. The set of all prime numbers which do not belong to $\pi$ will be denoted by $\pi'$. A positive integer $n$ will be called a {\color{blau} $\pi$-{\em number}} if each prime divisor of $n$ belongs to $\pi$. \na

\ne Let $S$ be an association scheme. An element of $S$ is called {\color{blau} $\pi$-{\em valenced}} if its valency is a $\pi$-number. A closed subset of $S$ is called {\color{blau} $\pi$-{\em valenced}} if each of its elements is $\pi$-valenced. A $\pi$-valenced closed subset $T$ of $S$ is called a closed {\color{blau} $\pi$-{\em subset}} if its valency is a $\pi$-number. A closed $\pi$-subset of $S$ is called a {\color{blau} \em Hall $\pi$-subset} of $S$ if its index in $S$ is a $\pi'$-number; cf.\ {\color{gruen} [8; Section 4.4]}. \na

\ne Let $T$ and $U$ be closed subsets of a scheme $S$. We say that $T$ and $U$ are {\color{blau} \em conjugate in} $S$ if $S$ contains an element $s$ with $s^{*}Ts=U$. Assume that $T\subseteq U$. In {\color{gruen} [8; Lemma 2.3.6(ii)]}, it was shown that $n_{T}$ divides $n_{U}$. The quotient $n_{U}/n_{T}$ is called the {\color{blau} \em index} of $T$ {\color{blau} \em in} $U$. The closed subset $T$ is said to be {\color{blau} \em strongly normal} in $U$ if $u^{*}Tu\subseteq T$ for each element $u$ in $U$. \na

\ne An association scheme $S$ will be called {\color{blau} \em solvable} if it contains closed subsets $T_{0}$, $\ldots$, $T_{n}$ such that $T_{0}=\{1\}$, $T_{n}=S$, and, for each element $i$ in $\{1,\ldots,n\}$, $T_{i-1}$ is strongly normal in $T_{i}$ and has prime index in $T_{i}$. \na

\ne It is easy to see that, via the group correspondence, a thin association scheme corresponds to a solvable group if and only if it is solvable, so that the notion of a solvable association scheme which we suggest here naturally generalizes the notion of a finite solvable group. \na

\ne Now we are ready to state the main result of this note. \gra

\ns\ne {\bf Theorem} \na

\ne {\it Let $\pi$ be a set of prime numbers, and let $S$ be a solvable and $\pi$-valenced association scheme. Then we have the following.} \tabulaturna

\ne\loss

\hfill (i) & {\it The scheme $S$ possesses at least one Hall $\pi$-subset.} \\

\end{tabular} \tabulaturkla

\ne\loss

\hfill (ii) & {\it Any two Hall $\pi$-subsets of $S$ are conjugate in $S$.} \\

\end{tabular} \tabulaturkla

\ne\loss

\hfill (iii) & {\it Any closed $\pi$-subset of $S$ is contained in a Hall $\pi$-subset of $S$.} \\

\end{tabular} \gra

\ne Via the group correspondence, the restriction of our theorem to the class of all thin association schemes results precisely in Hall's theorem {\color{gruen} [2; Theorem]} on finite solvable groups, since thin association schemes are $\pi$-valenced for any set $\pi$ of prime numbers. \na

\ns\ne It is perhaps not without interest to point out that both conditions under which we prove our theorem are needed, the solvability of $S$ as well as the hypothesis that $S$ is $\pi$-valenced. As for the necessity of the solvability we refer of course to Hall's note {\color{gruen} [2]}. The necessity of the second hypothesis is shown by the association scheme $HM_{176}(28)$ in the list {\color{gruen} [3]} of Akihide Hanaki and Izumi Miyamoto. This scheme is solvable and $\{2\}$-valenced, and it has valency $28$. It has a closed subset of valency $4$, but no closed subset of valency $7$. \na

\ns\ne The key to the proof of our main result is {\color{rot} Theorem 8.5}. In this theorem, we show that, given a set $\pi$ of prime numbers, each solvable $\pi$-valenced association scheme $S$ contains a strongly normal closed $\pi$-subset $U$. Since $U$ is strongly normal in $S$, we obtain from {\color{gruen} [8; Lemma 4.2.5(ii)]} that the quotient scheme $S/\!/U$ is thin, and our main result is then obtained by an application of Hall's above-mentioned theorem to the group which, via the group correspondence, corresponds to $S/\!/U$. \na

\ne The proof of {\color{rot} Theorem 8.5} is quite involved and relies on a series of results on hypergroups. We define hypergroups in the beginning of {\color{rot} Section 2}. The reason why our hypergroup theoretic efforts (in {\color{rot} Section 2} to {\color{rot} Section 7}) lead to a proof of {\color{rot} Theorem 8.5} is the observation that each association scheme $S$ is a hypergroup with respect to the hypermultiplication which assigns any pair $(p,q)$ of elements of $S$ to its complex product $pq$; cf.\ {\color{gruen} [8; Lemma 1.3.1]}, {\color{gruen} [8; Lemma 1.3.3(ii)]}, and {\color{gruen} [8; Lemma 1.3.3(i)]}. The set $S$ endowed with the above hypermultiplication will be called the {\color{blau} \em hypergroup defined by the complex multiplication in} $S$. \gra\gra

\ns\centerline{\bf 2. The definition of a hypergroup and some preliminaries on hypergroups} \gra\na

\ne We define a {\color{blau} \em hypermultiplication} on a set $S$ to be a map from $S\times S$ to the power set of $S$. The image of a pair $(p,q)$ of elements of a set $S$ under a hypermultiplication on $S$ will be denoted by $pq$, and it will be called the {\color{blau} \em hyperproduct} or simply the {\color{blau} \em product} of $p$ and $q$. \na

\ne Let $S$ be a set endowed with a hypermultiplication, and let $P$ and $Q$ be subsets of $S$. We write $PQ$ for the union of the products $pq$ with $p\in P$ and $q\in Q$, and we refer to $PQ$ as the {\color{blau} \em hyperproduct} or simply as the {\color{blau} \em product} of $P$ and $Q$. The map which takes any pair $(P,Q)$ of subsets of $S$ to the product $PQ$ will be called the {\color{blau} \em extension} of the hypermultiplication on $S$. \na

\ne If $P$ contains an element $p$ with $P=\{p\}$, we write $pQ$ instead of $PQ$. Similarly, if $Q$ contains an element $q$ with $Q=\{q\}$, we write $Pq$ instead of $PQ$. Note that $pq=PQ$ if $P$ contains an element $p$ with $P=\{p\}$ and $Q$ contains an element $q$ with $Q=\{q\}$. \na

\ns\ne Following (and slightly generalizing) Fr$\acute{{\rm e}}$d$\acute{{\rm e}}$ric Marty's terminology in {\color{gruen} [5]} we call a set $S$ endowed with a hypermultiplication a {\color{blau} \em hypergroup} if it satisfies the following conditions. \gra

\ne\loss

\ \ ${\rm H}1$ & For any three elements $p$, $q$, and $r$ in $S$, we have $p(qr)=(pq)r$. \\

\end{tabular} \tabulaturkla

\ne\loss

\ \ ${\rm H}2$ & The set $S$ contains an element $e$ such that $se=\{s\}$ for each element $s$ in $S$. \\

\end{tabular} \tabulaturkla

\ne\loss

\ \ ${\rm H}3$ & For each element $s$ in $S$, there exists an element $s^{*}$ in $S$ such that $q\in p^{*}r$ and $p\in rq^{*}$ for any three elements $p$, $q$, and $r$ in $S$ satisfying $r\in pq$. \\

\end{tabular} \na\na

\ne In {\color{rot} Lemma 2.3}, we will see that in a hypergroup $H$, the set $ab$ is not empty for any two elements $a$ and $b$ in $H$. \na

\ne An element $e$ of a hypergroup $H$ which satisfies $he=\{h\}$ for each element $h$ in $H$ is called a {\color{blau} \em neutral element} of $H$. It is easy to see and was shown in {\color{gruen} [9; Lemma 2.4]} that hypergroups cannot have two distinct neutral elements. This allows us to speak about {\em the} neutral element of a hypergroup. Within this article, the neutral element of a hypergroup will always be denoted by $1$. \na

\ns\ne A map $^{*}$ from a hypergroup $H$ to itself which satisfies $b\in a^{*}c$ and $a\in cb^{*}$ for any three elements $a$, $b$, and $c$ in $H$ with $c\in ab$ is called an {\color{blau} \em inverse function} of $H$. It is easy to see and was shown in {\color{gruen} [9; Lemma 2.7]} that hypergroups cannot have two different inverse functions. This allows us to speak about {\em the} inverse function of a hypergroup $H$. Within this article, the inverse of a hypergroup element $h$ will always be denoted by $h^{*}$. \gra

\ns\ne {\bf Lemma 2.1} \na

\ne {\it Let $H$ be a hypergroup, and let $a$ and $b$ be elements of $H$. Then $1\in a^{*}b$ if and only if $a=b$.} \gra

\ne {\it Proof.} Assume first that $1\in a^{*}b$. Then, by {\color{rot} Condition ${\rm H}3$}, $b\in a\!\cdot\! 1$. On the other hand, by {\color{rot} Condition ${\rm H}2$}, $a\!\cdot\! 1=\{a\}$. Thus, $a=b$. \na

\ne Assume that $a=b$. Then, by {\color{rot} Condition ${\rm H}2$}, $b\in a\!\cdot\! 1$. Thus, by {\color{rot} Condition ${\rm H}3$}, $1\in a^{*}b$. \qed\gra

\ns\ne {\bf Lemma 2.2} \na

\ne {\it Let $H$ be a hypergroup, and let $h$ be an element in $H$. Then $h^{**}=h$.} \gra

\ne {\it Proof.} From {\color{rot} Lemma 2.1} we know that $1\in h^{*}h$. It follows that $h\in h^{**}\!\cdot\! 1$; cf.\ {\color{rot} Condition ${\rm H}3$}. Now recall that, by {\color{rot} Condition ${\rm H}2$}, $h^{**}\!\cdot\! 1=\{h^{**}\}$. Thus, $h\in\{h^{**}\}$, and that implies that $h^{**}=h$. \qed\gra

\ns\ne The following lemma was communicated to the second author by Chris French. \gra

\ns\ne {\bf Lemma 2.3} \na

\ne {\it For any two elements $a$ and $b$ in $H$, the set $ab$ is not empty.} \gra

\ne {\it Proof.} Let $a$ and $b$ be elements in $H$. From Condition ${\rm H}2$, Lemma 2.1, and Condition ${\rm H}1$ we obtain that
\[
a\in a\!\cdot\! 1\subseteq a(b^{**}b^{*})=a(bb^{*})=(ab)b^{*}.
\]
This shows that $ab$ is not empty. \qed\gra

\ns\ne {\bf Lemma 2.4} \na

\ne {\it Let $H$ be a hypergroup, and let $a$, $b$, and $c$ be elements of $H$. Then the statements $c\in ab$, $b\in a^{*}c$, $a^{*}\in bc^{*}$, $c^{*}\in b^{*}a^{*}$, $b^{*}\in c^{*}a$, and $a\in cb^{*}$ are pairwise equivalent.} \gra

\ne {\it Proof.} From $c\in ab$ one obtains that $b\in a^{*}c$; cf.\ {\color{rot} Condition ${\rm H}3$}. Similarly, one obtains $a^{*}\in bc^{*}$ from $b\in a^{*}c$, $c^{*}\in b^{*}a^{*}$ from $a^{*}\in bc^{*}$, $b^{*}\in c^{*}a$ from $c^{*}\in b^{*}a^{*}$, $a\in cb^{*}$ from $b^{*}\in c^{*}a$, and $c\in ab$ from $a\in cb^{*}$. \qed\gra

\ns\ne Recall from {\color{rot} Lemma 2.1} that $1\in h^{*}h$ for each hypergroup element $h$. A hypergroup element $h$ is called {\color{blau} \em thin} if $h^{*}h=\{1\}$, and subsets of hypergroups will be called {\color{blau} \em thin} if all of their elements are thin. \gra

\ns\ne {\bf Lemma 2.5} \na

\ne {\it Let $H$ be a hypergroup, and let $a$ and $b$ be elements in $H$. Assume that $b$ is thin. Then $|ab|=1$.} \gra

\ne {\it Proof.} Let $c$ be an element in $ab$. Then $a\in cb^{*}$. Thus, as $b$ is thin, $ab\subseteq cb^{*}b=\{c\}$. \qee\gra

\ns\ne For each subset $A$ of a hypergroup, we set $A^{*}:=\{a^{*}\st a\in A\}$, and we notice that, by {\color{rot} Lemma 2.2}, $A^{**}=A$ for each subset $A$ of a hypergroup. \gra

\ns\ne {\bf Lemma 2.6} \na

\ne {\it Let $H$ be a hypergroup, and let $A$ and $B$ be subsets of $H$. Then the following hold.} \tabulaturna

\ne\loss

\hfill (i) & {\it Assume that $A\subseteq B$. Then $A^{*}\subseteq B^{*}$.} \\

\end{tabular} \tabulaturkla

\ne\loss

\hfill (ii) & {\it We have $(AB)^{*}=B^{*}A^{*}$.} \\

\end{tabular} \gra

\ne {\it Proof.} (i) Let $h$ be an element in $A^{*}$. Then $A$ contains an element $a$ with $h=a^{*}$. Since $a\in A$, $a\in B$. Thus, $h=a^{*}\in B^{*}$. \na

\ne (ii) Let $h$ be an element in $(AB)^{*}$. Then $AB$ contains an element $c$ with $h=c^{*}$. Since $c\in AB$, there exist elements $a$ in $A$ and $b$ in $B$ such that $c\in ab$. Thus, by {\color{rot} Lemma 2.4}, $c^{*}\in b^{*}a^{*}\subseteq B^{*}A^{*}$. Since $h=c^{*}$ this implies that $h\in B^{*}A^{*}$. \qed\gra\gra

\ns\centerline{\bf 3. Closed subsets of hypergroups} \gra\na

\ne A non-empty subset $A$ of a hypergroup is called {\color{blau} \em closed}\index{closed subset} if $A^{*}A\subseteq A$. \na

\ne Note that a closed subset $F$ of a hypergroup $H$ is a hypergroup with respect to the hypermultiplication which one obtains from the hypermultiplication of $H$ if one restricts the domain of the hypermultiplication of $H$ to $F\times F$ and the codomain of the hypermultiplication of $H$ to the power set of $F$. \gra

\ns\ne {\bf Lemma 3.1} \na

\ne {\it A subset $A$ of a hypergroup is closed if and only if $1\in A$, $A^{*}=A$, and $AA=A$.} \gra

\ne {\it Proof.} Assume first that $A$ is closed. Then, by definition, $A$ is not empty. Let $a$ be an element in $A$. Then, by definition, $a^{*}a\subseteq A$. On the other hand, by {\color{rot} Lemma 2.1}, $1\in a^{*}a$. Thus, $1\in A$. \na

\ne From $1\in A$ we obtain that $a^{*}\in a^{*}\!\cdot\! 1\subseteq A^{*}A\subseteq A$ for each element $a$ in $A$. This shows that $A^{*}\subseteq A$. From this it follows that $A=A^{**}\subseteq A^{*}$; cf.\ {\color{rot} Lemma 2.6(i)}. It follows that $A^{*}=A$. \na

\ne From $1\in A$ we obtain that $A=A\!\cdot\! 1\subseteq AA$. From $A^{*}=A$ we obtain that $AA=A^{*}A\subseteq A$. Thus, $AA=A$. \na

\ne Assume now, conversely, that $1\in A$, $A^{*}=A$, and $AA=A$. From $1\in A$ we obtain that $A$ is not empty. From $A^{*}\subseteq A$ and $AA\subseteq A$ we obtain that $A^{*}A\subseteq A$. \qed\gra

\ns\ne {\bf Lemma 3.2} \na

\ne {\it Let $H$ be a hypergroup, and let ${\cal F}$ be a non-empty set of closed subsets of $H$. Then the intersection of the closed subsets which belong to ${\cal F}$ is a closed subset of $H$.} \gra

\ne {\it Proof.} Let $A$ denote the intersection of the closed subsets of $H$ which belong to ${\mathcal F}$. From {\color{rot} Lemma 3.1} we know that $1\in F$ for each element $F$ in ${\cal F}$. Thus, $1\in A$, and that shows that $A$ is not empty. \na

\ne Let $F$ be an element in ${\cal F}$. Then $A\subseteq F$. Thus, by {\color{rot} Lemma 2.6(i)}, $A^{*}\subseteq F^{*}$. Thus, as $A\subseteq F$, we obtain that $A^{*}A\subseteq F^{*}F\subseteq F$, and since $F$ has been chosen arbitrarily from ${\cal F}$, this proves that $A^{*}A\subseteq A$. \qed\gra

\ns\ne {\bf Lemma 3.3} \na

\ne {\it Let $H$ be a hypergroup, let $F$ be a closed subset of $H$, and let $A$ and $B$ be subsets of $H$. Then we have the following.} \tabulaturna

\ne\loss

\hfill (i) & {\it If $A\subseteq F$, $A(B\cap F)=AB\cap F$.} \\

\end{tabular} \tabulaturkla

\ne\loss

\hfill (ii) & {\it If $B\subseteq F$, $(A\cap F)B=AB\cap F$.} \\

\end{tabular} \gra

\ne {\it Proof.} (i) Assume that $A\subseteq F$. Then, by {\color{rot} Lemma 3.1}, $A(B\cap F)\subseteq F$. Thus, as $A(B\cap F)\subseteq AB$, $A(B\cap F)\subseteq AB\cap F$. \na

\ne To show the reverse containment, let $f$ be an element in $AB\cap F$. Since $f\in AB$, there exist elements $a$ in $A$ and $b$ in $B$ such that $f\in ab$. From $f\in ab$ we obtain that $b^{*}\in f^{*}a$; cf.\ {\color{rot} Lemma 2.4}. \na

\ne Since $F$ is a closed subset of $H$, we obtain from $f\in F$ and $a\in A\subseteq F$ that $f^{*}a\subseteq F$. Thus, as $b^{*}\in f^{*}a$, $b^{*}\in F$. Since $F$ is closed, this implies that $b\in F$. Thus, as $b\in B$, $b\in B\cap F$. Now, as $f\in ab$, $f\in A(B\cap F)$. \na

\ne (ii) Assume that $B\subseteq F$. Then $B^{*}\subseteq F$. Thus, by {\color{rot} (i)}, $B^{*}(A^{*}\cap F)=B^{*}A^{*}\cap F$. Thus, as $F^{*}=F$, the claim follows from {\color{rot} Lemma 2.6(ii)}. \qed\gra

\ns\ne {\bf Lemma 3.4} \na

\ne {\it Let $H$ be a hypergroup, and let $D$ and $E$ be closed subsets of $H$. Then the following hold.} \tabulaturna

\ne\loss

\hfill (i) & {\it For any two elements $a$ and $b$ in $H$ with $a\in DbE$, we have $DbE\subseteq DaE$.} \\

\end{tabular} \tabulaturkla

\ne\loss

\hfill (ii) & {\it The set $\{DhE\st h\in H\}$ is a partition of $H$.} \\

\end{tabular} \gra

\ne {\it Proof.} (i) Let $a$ and $b$ be elements in $H$, and assume that $a\in DbE$. Then there exist elements $e$ in $D$ and $f$ in $E$ such that $a\in ebf$. From $a\in ebf$ we obtain an element $c$ in $eb$ such that $a\in cf$. From $c\in eb$ we obtain that $b\in e^{*}c$, from $a\in cf$ we obtain that $c\in af^{*}$. Thus, $b\in e^{*}af^{*}\subseteq DaE$. Since $D$ and $E$ are assumed to be closed, this implies that $DbE\subseteq DaE$; cf.\ {\color{rot} Lemma 3.1}. \na

\ns\ne (ii) Since $1\in D$ and $1\in E$, we have $h\in DhE$ for each element $h$ in $H$. Thus, $H$ is equal to the union of the sets $DhE$ with $h\in H$. \na

\ne To show that $\{DhE\st h\in H\}$ is a partition of $H$ we now choose elements $a$ and $b$ in $H$, and we assume that $DaE\cap DbE$ is not empty. We have to show that $DaE=DbE$. \na

\ne Since $DaE\cap DbE$ is not empty, we find an element $c$ in $DaE\cap DbE$. Since $c\in DaE$, $DcE\subseteq DaE$. From $c\in DaE$ we also obtain that $DaE\subseteq DcE$; cf.\ {\color{rot} (i)}. Thus, $DaE=DcE$. \na

\ne Similarly, $DbE=DcE$, so that $DaE=DbE$. \qed\gra

\ns\ne {\bf Lemma 3.5} \na

\ne {\it Let $H$ be a hypergroup, and let $D$ and $E$ be closed subsets of $H$. Then $DE$ is a closed subset of $H$ if and only if $DE=ED$.} \gra

\ne {\it Proof.} Assume first that $DE$ is closed. Then, by {\color{rot} Lemma 3.1}, $(DE)^{*}=DE$. On the other hand, as $D$ and $E$ are closed, we also have $D^{*}=D$ and $E^{*}=E$. Thus, by {\color{rot} Lemma 2.6(ii)}, $(DE)^{*}=ED$. From $(DE)^{*}=DE$ and $(DE)^{*}=ED$ we obtain that $DE=ED$. \na

\ne Assume, conversely, that $DE=ED$. Then referring to {\color{rot} Lemma 2.6(ii)} we obtain that
\[
(DE)^{*}DE=E^{*}D^{*}DE\subseteq E^{*}DE=E^{*}ED\subseteq ED=DE.
\]
Thus, as $DE$ is not empty, $DE$ is closed. \qed\gra\gra

\ns\centerline{\bf 4. Normality} \gra\na

\ne Let $H$ be a hypergroup, and let $D$ and $E$ be closed subsets of $H$. We say that $D$ {\color{blau} \em normalizes} $E$ if $Ed\subseteq dE$ for each element $d$ in $D$. If $H$ normalizes a closed subset $F$ of $H$, we say that $F$ is {\color{blau} \em normal in} $H$. \gra

\ns\ne {\bf Lemma 4.1} \na

\ne {\it Let $H$ be a hypergroup, and let $D$ and $E$ be closed subsets of $H$. Assume that $D$ normalizes $E$. Then the following hold.} \tabulaturna

\ne\loss

\hfill (i) & {\it For each element $d$ in $D$, we have $Ed=dE$.} \\

\end{tabular} \tabulaturkla

\ne\loss

\hfill (ii) & {\it The product $ED$ is a closed subset of $H$.} \\

\end{tabular} \tabulaturkla

\ne\loss

\hfill (iii) & {\it The closed subset $E$ of $H$ is normal in $ED$.} \\

\end{tabular} \tabulaturkla

\ne\loss

\hfill (iv) & {\it The closed subset $E\cap D$ is normal in $D$.} \\

\end{tabular} \gra

\ne {\it Proof.} (i) Let $d$ be an element in $D$. Then, $d^{*}\in D$. Thus, as $D$ is assumed to normalize $E$, $Ed\subseteq dE$ and $Ed^{*}\subseteq d^{*}E$. Since $E$ is closed, $Ed^{*}\subseteq d^{*}E$ implies that $dE\subseteq Ed$; cf.\ {\color{rot} Lemma 2.6}. \na

\ne (ii) From {\color{rot} (i)} we obtain that $ED=DE$. Thus, by {\color{rot} Lemma 3.5}, $ED$ is closed. \na

\ns\ne (iii) Let $h$ be an element in $ED$. Then $D$ contains an element $d$ with $h\in Ed$. Since $D$ normalizes $E$, we also have $Ed\subseteq dE$. Thus, $h\in dE$. Now, by {\color{rot} Lemma 3.4(ii)}, $hE=dE$. Thus,
\[
Eh\subseteq Ed\subseteq dE=hE,
\]
and that means that $h$ normalizes $E$. Since $h$ has been chosen arbitrarily in $ED$, we have seen that $ED$ normalizes $E$. \na

\ne (iv) We are assuming that $D$ normalizes $E$. Thus, by definition, $Ed\subseteq dE$ for each element $d$ in $D$. It follows that
\[
(E\cap D)d=Ed\cap D\subseteq dE\cap D=d(E\cap D)
\]
for each element $d$ in $D$; cf.\ {\color{rot} Lemma 3.3}. \qed\gra

\ns\ne A closed subset $E$ of a hypergroup $H$ is said to be {\color{blau} \em subnormal in} $H$ if $H$ contains closed subsets $F_{0}$, $\ldots$, $F_{n}$ such that $F_{0}=E$, $F_{n}=H$, and, for each element $i$ in $\{1,\ldots,n\}$, $F_{i-1}$ is normal in $F_{i}$. \gra

\ns\ne {\bf Lemma 4.2} \na

\ne {\it Let $H$ be a hypergroup, and let $D$ and $E$ be closed subsets of $H$. Assume that $D$ is subnormal in $H$ and that $E$ is normal in $H$. Then $ED$ is a closed subset of $H$ and is subnormal in $H$.} \gra

\ne {\it Proof.} That $ED$ is a closed subset of $H$ follows from {\color{rot} Lemma 4.1(ii)}. \na

\ne We are assuming that $D$ is subnormal in $H$. Thus, there exist closed subsets $F_{0}$, $\ldots$, $F_{n}$ of $H$ such that $F_{0}=D$, $F_{n}=H$ and, for each element $i$ in $\{1,\ldots,n\}$, $F_{i-1}$ is normal in $F_{i}$. We will see that, for each element $i$ in $\{1,\ldots,n\}$, $EF_{i-1}$ is normal in $EF_{i}$. \na

\ne Let $h$ be an element in $EF_{i}$. Then $F_{i}$ contains an element $f$ such that $h\in Ef$. It follows that $Eh\subseteq Ef$. From $h\in Ef$ and $Ef\subseteq fE$ we also obtain that $h\in fE$, so that, by {\color{rot} Lemma 3.4(ii)}, $hE=fE$. It follows that
\[
EF_{i-1}h=F_{i-1}Eh\subseteq F_{i-1}Ef\subseteq F_{i-1}fE\subseteq fF_{i-1}E=fEF_{i-1}=hEF_{i-1}.
\]
Since $h$ has been chosen arbitrarily from $EF_{i}$, this shows that $EF_{i-1}$ is normal in $EF_{i}$. \qed\gra

\ns\ne A closed subset $F$ of a hypergroup $H$ is called {\color{blau} \em strongly normal in} $H$ if $h^{*}Fh\subseteq F$ for each element $h$ in $H$. \gra

\ns\ne {\bf Lemma 4.3} \na

\ne {\it Let $H$ be a hypergroup, and let $F$ be a closed subset of $H$. Assume that $F$ is strongly normal in $H$. Then $F$ is normal in $H$.} \gra

\ne {\it Proof.} Since $F$ is assumed to be strongly normal in $H$, we have $h^{*}Fh\subseteq F$ for each element $h$ in $H$. It follows that $Fh\subseteq hh^{*}Fh\subseteq hF$ for each element $h$ in $H$, and that means that $F$ is normal in $H$. \qed\gra

\ns\ne For each hypergroup $H$, we denote by $O^{\vartheta}(H)$ the intersection of all closed subsets of $H$ which are strongly normal in $H$. \gra

\ns\ne {\bf Lemma 4.4} \na

\ne {\it Let $H$ be a hypergroup. Then the following hold.} \tabulaturna

\ne\loss

\hfill (i) & {\it The set $O^{\vartheta}(H)$ is a closed subset of $H$ and strongly normal in $H$.} \\

\end{tabular} \tabulaturkla

\ne\loss

\hfill (ii) & {\it For each element $h$ in $H$, we have $h^{*}h\subseteq O^{\vartheta}(H)$.} \\

\end{tabular} \gra

\ne {\it Proof.} (i) From {\color{rot} Lemma 3.2} we know that $O^{\vartheta}(H)$ is a closed subset of $H$. \na

\ne Let ${\cal F}$ denote the set of all closed subsets of $H$ which are strongly normal in $H$, and let $F$ be an element in ${\cal F}$. Then we have
\[
h^{*}O^{\vartheta}(H)h\subseteq h^{*}Fh\subseteq F
\]
for each element $h$ in $H$. Thus, $h^{*}O^{\vartheta}(H)h\subseteq O^{\vartheta}(H)$ for each element $h$ in $H$, and that means that $O^{\vartheta}(H)$ is strongly normal in $H$. \na

\ne (ii) For each element $h$ in $H$, we have $h^{*}h\subseteq h^{*}O^{\vartheta}(H)h\subseteq O^{\vartheta}(H)$; cf.\ {\color{rot} (i)}. \qed\gra

\ns\ne The set of all thin elements of a hypergroup $H$ will be denoted by $O_{\vartheta}(H)$, and a hypergroup $H$ will be called {\color{blau} \em metathin} if $O^{\vartheta}(H)\subseteq O_{\vartheta}(H)$. \gra

\ns\ne {\bf Lemma 4.5} \na

\ne {\it Let $H$ be a metathin hypergroup, and let $h$ be an element in $H$. Then the following hold.} \tabulaturna

\ne\loss

\hfill (i) & {\it We have $\{h\}=hh^{*}h$.} \\

\end{tabular} \tabulaturkla

\ne\loss

\hfill (ii) & {\it The set $h^{*}h$ is a closed subset of $H$.} \\

\end{tabular} \tabulaturkla

\ne\loss

\hfill (iii) & {\it The set $h^{*}h$ is thin.} \\

\end{tabular} \tabulaturkla

\ne\loss

\hfill (iv) & {\it The set $h^{*}h$ is normal in $O^{\vartheta}(H)$.} \\

\end{tabular} \gra

\ne {\it Proof.} (i) Let $a$ be an element in $hh^{*}h$. Then $h^{*}h$ contains an element $b$ such that $a\in hb$. Since $b\in h^{*}h$, $b\in O^{\vartheta}(H)$; cf.\ {\color{rot} Lemma 4.4(ii)}. Thus, as $H$ is assumed to be metathin, $b$ is thin, and that means that $b^{*}b=\{1\}$. \na

\ne From $a\in hb$ we obtain that $h\in ab^{*}$, from $b\in h^{*}h$ we obtain that $h^{*}\in bh^{*}$. Thus, as $b^{*}b=\{1\}$, we obtain from {\color{rot} Lemma 2.1} that $1\in hh^{*}\subseteq ab^{*}bh^{*}=ah^{*}$. Thus, by {\color{rot} Lemma 2.1}, $a^{*}=h^{*}$, whence $a=h$. \na

\ne Since $a$ has been chosen arbitrarily from $hh^{*}h$, we have shown that $hh^{*}h\subseteq\{h\}$. On the other hand, $hh^{*}h$ is not empty. Thus, $\{h\}=hh^{*}h$. \na

\ns\ne (ii) From {\color{rot} Lemma 2.6(ii)} we know that $(h^{*}h)^{*}=h^{*}h$. Thus, by {\color{rot} (i)}, $(h^{*}h)^{*}h^{*}h=h^{*}hh^{*}h=h^{*}h$. Thus, as $h^{*}h$ is not empty, $h^{*}h$ is a closed subset of $H$. \na

\ne (iii) We have $h^{*}h\subseteq O^{\vartheta}(H)\subseteq O_{\vartheta}(H)$. Thus, $h^{*}h$ is thin. \na

\ne (iv) In order to show that $h^{*}h$ is normal in $O^{\vartheta}(H)$, we choose an element $a$ in $O^{\vartheta}(H)$. We have to show that $a^{*}h^{*}ha\subseteq h^{*}h$. \na

\ne Since $H$ is assumed to be metathin, we have $O^{\vartheta}(H)\subseteq O_{\vartheta}(H)$. Thus, as $a\in O^{\vartheta}(H)$, $a$ is thin. Thus, by {\color{rot} Lemma 2.5}, $|ha|=1$. Now recall from {\color{rot} Lemma 4.4(i)} that $O^{\vartheta}(H)$ is strongly normal in $H$. Thus, $O^{\vartheta}(H)$ contains an element $b$ such that $ha\subseteq bh$. Since $b\in O^{\vartheta}(H)$ and $O^{\vartheta}(H)\subseteq O_{\vartheta}(H)$, $b$ is thin, so $b^{*}b=\{1\}$. It follows that $a^{*}h^{*}ha\subseteq h^{*}b^{*}bh=h^{*}h$; cf.\ {\color{rot} Lemma 2.6}. \qed\gra\gra

\ns\centerline{\bf 5. Quotients of hypergroups} \gra\na

\ne Let $H$ be a hypergroup, and let $F$ be a closed subset of $H$. \na

\ne For each element $h$ of $H$, we set
\[
h^{F}:=FhF,
\]
and we define $H/\!/F$ to be the set of all hyperproducts $h^{F}$ with $h\in H$. In {\color{gruen} [7; Lemma 2.2]}, it was shown that we obtain a hypermultiplication on $H/\!/F$ if we define
\[
a^{F}b^{F}:=\{h^{F}\st h\in aFb\}
\]
for any two elements $a$ and $b$ in $H$. We call this hypermultiplication on $H/\!/F$ the {\color{blau} \em hypermultiplication} on $H/\!/F$ {\color{blau} \em defined by} $F$. \na

\ns\ne We notice that the hypermultiplication on $H/\!/F$ defined by $F$ carries the danger of ambiguity. In fact, the hypermultiplication on $H/\!/F$ defined by $F$ is defined on the same set as the restriction to $H/\!/F\times H/\!/F$ of the extension of the hypermultiplication of $H$. However, while the codomain of the restriction to $H/\!/F\times H/\!/F$ of the extension of the hypermultiplication of $H$ is the power set of $H/\!/F$, the codomain of the hypermultiplication on $H/\!/F$ defined by $F$ is the power set of $H$. To say it differently, for any two elements $a$ and $b$ in $H$, $(FaF)(FbF)$ stands for the set union of the products $FhF$ with $h\in aFb$, whereas $a^{F}b^{F}$ is defined to be the of all products $h^{F}$ with $h\in aFb$, although $FaF=a^{F}$ and $FbF=b^{F}$. In the following, we will take care that no misunderstanding will arise in this regard. \na

\ne In {\color{gruen} [7; Lemma 2.3]}, it was shown that $H/\!/F$ is a hypergroup with respect to the hypermultiplication on $H/\!/F$ defined by $F$ (and with inverse function $h^{F}\mapsto (h^{*})^{F}$). We refer to $H/\!/F$ as the {\color{blau} \em quotient} of $H$ {\color{blau} \em over} $F$. \gra

\ns\ne {\bf Lemma 5.1} \na

\ne {\it Let $H$ be a hypergroup, let $F$ be a closed subset of $H$, let $h$ be an element in $H$, and let $A_{1}$, $\ldots$, $A_{n}$ be subsets of $H$. Then we have $h^{F}\in (A_{1}/\!/F)\cdots (A_{n}/\!/F)$ if and only if $h\in (FA_{1}F)\cdots (FA_{n}F)$.} \gra

\ne {\it Proof.} We first prove the statement for $n=1$. By definition, we have $h^{F}\in A_{1}/\!/F$ if and only if there exists an element $a$ in $A_{1}$ such that $h^{F}=a^{F}$. By definition, we also have $h^{F}=a^{F}$ if and only if $FhF=FaF$. Thus, $h^{F}\in A_{1}/\!/F$ if and only if $h\in FA_{1}F$; cf.\ {\color{rot} Lemma 3.4(ii)}. \na

\ne Now we assume that $2\le n$, and we first suppose that $h^{F}\in (A_{1}/\!/F)\cdots (A_{n}/\!/F)$. Then, by definition, there exist elements $b$ and $c$ in $H$ such that $b^{F}\in (A_{1}/\!/F)\cdots (A_{n-1}/\!/F)$, $c^{F}\in A_{n}/\!/F$, and $h^{F}\in b^{F}c^{F}$. From $h^{F}\in b^{F}c^{F}$ we obtain an element $d$ in $bFc$ such that $h^{F}=d^{F}$. It follows that $h\in FdF\subseteq FbFcF$. \na

\ns\ne Since $b^{F}\in (A_{1}/\!/F)\cdots (A_{n-1}/\!/F)$, induction yields $b\in (FA_{1}F)\cdots (FA_{n-1}F)$. In the first paragraph of this proof, we also saw that $c^{F}\in A_{n}/\!/F$ implies that $c\in FA_{n}F$. Thus, we have
\[
FbFcF\subseteq (FA_{1}F)\cdots (FA_{n}F).
\]
Since $h\in FbFcF$, this implies that $h\in (FA_{1}F)\cdots (FA_{n}F)$. \na

\ne Suppose, conversely, that $h\in (FA_{1}F)\cdots (FA_{n}F)$. Then, by definition,  there exist elements $b$ in $(FA_{1}F)\cdots (FA_{n-1}F)$ and $c$ in $FA_{n}F$ such that $h\in bc$. From $h\in bc$ (together with $bc\subseteq bFc$) we obtain that $h^{F}\in b^{F}c^{F}$. \na

\ne Since $b\in (FA_{1}F)\cdots (FA_{n-1}F)$, induction yields $b^{F}\in (A_{1}/\!/F)\cdots (A_{n-1}/\!/F)$. In the first part of this proof, we also saw that $c\in FA_{n}F$ implies that $c^{F}\in A_{n}/\!/F$. Thus, we have
\[
b^{F}c^{F}\subseteq (A_{1}/\!/F)\cdots (A_{n}/\!/F).
\]
Since $h^{F}\in b^{F}c^{F}$, this implies that $h^{F}\in (A_{1}/\!/F)\cdots (A_{n}/\!/F)$. \qed\gra

\ns\ne {\bf Theorem 5.2} \na

\ne {\it Let $H$ be a hypergroup, and let $D$ be a closed subset of $H$. Then $E\mapsto E/\!/D$ is a bijective map from the set of all closed subsets of $H$ containing $D$ to the set of all closed subsets of $H/\!/D$.} \gra

\ne {\it Proof.} Let $E$ be a closed subset of $H$, and assume that $D\subseteq E$. We first show that $E/\!/D$ is a closed subset of $H/\!/D$. \na

\ne Since $E$ is closed, we have $1\in E$; cf.\ {\color{rot} Lemma 3.1}. Thus, by definition, $1^{D}\in E/\!/D$. \na

\ne Since $E$ is closed, we also have $e^{*}\in E$ for each each element $e$ in $E$; cf.\ {\color{rot} Lemma 3.1}. Thus, by definition, $(e^{*})^{D}\in E/\!/D$. Since $(e^{D})^{*}=(e^{*})^{D}$, this implies that $(e^{D})^{*}\in E/\!/D$. \na

\ne Let $a$ and $b$ be elements in $E$, and let $h$ be an element in $H$ with $h^{D}\in a^{D}b^{D}$. From $h^{D}\in a^{D}b^{D}$ we obtain that $h\in DaDbD$; cf.\ {\color{rot} Lemma 5.1}. Thus, as $DaDbD\subseteq E$, $h\in E$, so that $h^{D}\in E/\!/D$. Since $h$ has been chosen arbitrarily from $H$ with $h^{D}\in a^{D}b^{D}$ , this shows that $a^{D}b^{D}\subseteq E/\!/D$. \na

\ns\ne What we have seen so far is that $E/\!/D$ is closed; cf.\ {\color{rot} Lemma 3.1}. \na

\ne Since $E$ has been chosen arbitrarily among the closed subset of $H$ containing $D$, this shows that $E\mapsto E/\!/D$ is a map from the set of all closed subsets of $H$ containing $D$ to the set of all closed subsets of $H/\!/D$. \na

\ne Let $B$ and $C$ be closed subsets of $H$ with $D\subseteq B$, $D\subseteq C$, and $B/\!/D=C/\!/D$. Let $b$ be an element in $B$. Then $b^{D}\in C/\!/D$. Thus, $C$ contains an element $c$ with $b^{D}=c^{D}$. It follows that $b\in DcD\subseteq C$. We, thus, have shown that $B\subseteq C$. Similarly, one shows that $C\subseteq B$, so that we have $B=C$. \na

\ns\ne Let $A$ be a subset of $H$, and assume that $A/\!/D$ is a closed subset of $H/\!/D$. From {\color{rot} Lemma 5.1} one obtains that $A/\!/D=DAD/\!/D$. Thus, we shall be done if we succeed in showing that $DAD$ is a closed subset of $H$. \na

\ne Let $h$ be an element in $DAD$. Then, by {\color{rot} Lemma 5.1}, $h^{D}\in A/\!/D$. Thus, as $A/\!/D$ is assumed to be closed, $(h^{D})^{*}\in A/\!/D$. Thus, as $(h^{D})^{*}=(h^{*})^{D}$, $(h^{*})^{D}\in A/\!/D$. Now, by definition, $(h^{*})^{D}=a^{D}$ for some element $a$ in $A$. It follows that $h^{*}\in DaD\subseteq DAD$. \na

\ne Let $b$ and $c$ be elements in $DAD$, and let $h$ be an element in $bc$. Then $h\in (DAD)(DAD)$. Thus, by {\color{rot} Lemma 5.1}, $h^{D}\in (A/\!/D)(A/\!/D)$. Since $A/\!/D$ is assumed to be closed, this implies that $h^{D}\in A/\!/D$. Thus, $A$ contains an element $a$ with $h^{D}=a^{D}$. It follows that $h\in DaD\subseteq DAD$. \qed\gra

\ns\ne {\bf Lemma 5.3} \na

\ne {\it Let $H$ be a hypergroup, let $D$ and $E$ be closed subsets of $H$, and assume that $D\subseteq E$. Then we have the following.} \tabulaturna

\ne\loss

\hfill (i) & {\it Assume that $E$ is normal in $H$. Then $E/\!/D$ is a normal closed subset of $H/\!/D$.} \\

\end{tabular} \tabulaturkla

\ne\loss

\hfill (ii) & {\it The closed subset $E$ of $H$ is strongly normal in $H$ if and only if $E/\!/D$ is strongly normal in $H/\!/D$.} \\

\end{tabular} \gra

\ne {\it Proof.} (i) Let $h$ be an element in $H$. Then, as $E$ is assumed to be normal in $H$, we have $Eh\subseteq hE$. It follows that $EhD\subseteq hED=hE\subseteq DhE$. Thus, by {\color{rot} Lemma 5.1}, $(E/\!/D)h^{D}\subseteq h^{D}(E/\!/D)$. \na

\ns\ne (ii) Let $h$ be an element in $H$. We first assume that $h^{*}Eh\subseteq E$ and will show that $(h^{D})^{*}(E/\!/D)h^{D}\subseteq E/\!/D$. \na

\ne Let $a$ be an element in $H$ with $a^{D}\in (h^{D})^{*}(E/\!/D)h^{D}$. Since $(h^{D})^{*}=(h^{*})^{D}$, we then have $a^{D}\in (h^{*})^{D}(E/\!/D)h^{D}$, so that, by {\color{rot} Lemma 5.1}, $a\in Dh^{*}EhD$. Thus, as $h^{*}Eh\subseteq E$ and $D\subseteq E$, we conclude that $a\in E$. It follows that $a^{D}\in E/\!/D$. \na

\ne Now we assume that $(h^{D})^{*}(E/\!/D)h^{D}\subseteq E/\!/D$ and will show that $h^{*}Eh\subseteq E$. \na

\ne Let $a$ be an element in $h^{*}Eh$. From $a\in h^{*}Eh$ (together with $(h^{*})^{D}=(h^{D})^{*}$) we obtain that $a^{D}\in (h^{D})^{*}(E/\!/D)h^{D}$; cf.\ {\color{rot} Lemma 5.1}. Thus, as $(h^{D})^{*}(E/\!/D)h^{D}\subseteq E/\!/D$, we have $a^{D}\in E/\!/D$. Since $D\subseteq E$, this implies that $a\in E$; cf.\ {\color{rot} Lemma 5.1}. \qed\gra

\ns\ne {\bf Lemma 5.4} \na

\ne {\it Let $H$ be a hypergroup, and let $F$ be a closed subset of $H$. Then $F$ is strongly normal in $H$ if and only if $H/\!/F$ is thin.} \gra

\ne {\it Proof.} Let $h$ be an element in $H$. Since $F$ is closed, we have $h^{*}Fh\subseteq F$ if and only if $(Fh^{*}F)(FhF)\subseteq F$. From {\color{rot} Lemma 5.1} we also know that $(Fh^{*}F)(FhF)\subseteq F$ if and only if $(h^{*})^{F}h^{F}\subseteq\{1^{F}\}$. Thus, as $(h^{*})^{F}=(h^{F})^{*}$, we also have that $(h^{*})^{F}h^{F}\subseteq\{1^{F}\}$ if and only if $h^{F}$ is thin. \qed\gra\gra

\ns\centerline{\bf 6. Homomorphisms of hypergroups} \gra\na

\ne Let $H$ and $H'$ be hypergroups. \na

\ne As it is customary, we set $\phi(A):=\{\phi(a)\st a\in A\}$ whenever $\phi$ is a map from $H$ to $H'$ and $A$ a subset of $H$. \na

\ne A map $\phi$ from $H$ to $H'$ is called a {\color{blau} \em homomorphism} if
\[
\phi(ab)=\phi(a)\phi(b)
\]
for any two elements $a$ and $b$ in $H$ and $\phi(1)=1$. \gra

\ns\ne {\bf Lemma 6.1} \na

\ne {\it Let $H$ and $H'$ be hypergroups, let $\phi$ be a homomorphism from $H$ to $H'$, and let $h$ be an element in $H$. Then $\phi(h^{*})=\phi(h)^{*}$.} \gra

\ne {\it Proof.} By definition, $\phi(1)=1$. On the other hand, by {\color{rot} Lemma 2.1}, $1\in h^{*}h$. Thus,
\[
1=\phi(1)\in\phi(h^{*}h)=\phi(h^{*})\phi(h),
\]
so that, by {\color{rot} Lemma 2.1}, $\phi(h^{*})=\phi(h)^{*}$. \qed\gra

\ns\ne Let $H$ and $H'$ be hypergroups, and let $\phi$ be a homomorphism from $H$ to $H'$. We define
\[
{\rm ker}(\phi):=\{h\in H\st\phi(h)=1\}
\]
and call this set the {\color{blau} \em kernel} of $\phi$. \gra

\ns\ne {\bf Lemma 6.2} \na

\ne {\it Let $H$ and $H'$ be hypergroups, let $\phi$ be a homomorphism from $H$ to $H'$, and set $F:={\rm ker}(\phi)$. Let $a$ and $b$ be elements in $F$. Then we have $\phi(a)=\phi(b)$ if and only if $aF=bF$.} \gra

\ne {\it Proof.} From {\color{rot} Lemma 6.1} we obtain that
\[
\phi(a)^{*}\phi(b)=\phi(a^{*})\phi(b)=\phi(a^{*}b).
\]
Thus, by {\color{rot} Lemma 2.1}, $\phi(a)=\phi(b)$ if and only if $1\in\phi(a^{*}b)$. \na

\ne On the other hand, we have $1\in\phi(a^{*}b)$ if and only if $a^{*}b$ contains an element $f$ with $\phi(f)=1$, and this means that $f\in F$. \na

\ne To conclude the proof we notice that $f\in a^{*}b$ if and only if $b\in af$; cf.\ {\color{rot} Lemma 2.2}. Thus, we have $1\in\phi(a^{*}b)$ if and only if $aF=bF$; cf.\ {\color{rot} Lemma 3.4(ii)}. \qed\gra

\ns\ne Bijective hypergroup homomorphisms are called {\color{blau} \em isomorphisms}. \na

\ne If there exists an isomorphism from $H$ to $H'$ (or from $H'$ to $H$), we say that $H$ and $H'$ are {\em isomorphic}\index{isomorphic}, and we indicate this by writing $H\cong H'$. \na

\ne Let $\phi$ be a homomorphism from $H$ to $H'$. As is customary, we define
\[
{\rm im}(\phi):=\phi(H)
\]
and call this set the {\color{blau} \em image} of $\phi$. \na

\ne The following two theorems are {\color{gruen} [7; Theorem 3.3]} and {\color{gruen} [7; Theorem 3.4(ii)]}. For the sake of completeness of this article, we include a proof. \gra

\ns\ne {\bf Theorem 6.3} \na

\ne {\it Let $H$ and $H'$ be hypergroups, and let $\phi$ be a homomorphism from $H$ to $H'$. Then ${\rm ker}(\phi)$ is a normal closed subset of $H$, ${\rm im}(\phi)$ is a closed subset of $H'$, and $H/\!/{\rm ker}(\phi)\cong {\rm im}(\phi)$.} \gra

\ne {\it Proof.} By definition, $\phi(1)=1$. Thus, $1\in {\rm ker}(\phi)$, so ${\rm ker}(\phi)$ is not empty. \na

\ne Let $a$ and $b$ be elements in ${\rm ker}(\phi)$. Then $\phi(a)=1$ and $\phi(b)=1$. Thus, referring to {\color{rot} Lemma 6.1} we obtain that
\[
\phi(a^{*}b)=\phi(a^{*})\phi(b)=\phi(a)^{*}\phi(b)=1^{*}1=\{1\},
\]
and that means that $a^{*}b\subseteq {\rm ker}(\phi)$. \na

\ns\ne So far, we have seen that ${\rm ker}(\phi)$ is a closed subset of $H$. In order to show that ${\rm ker}(\phi)$ is normal in $H$, we set $F:={\rm ker}(\phi)$. \na

\ne Let $f$ be an element in $F$, and let $h$ be an element in $H$. We have to show that $fh\subseteq hF$. \na

\ne Let $a$ be an element in $fh$. Then
\[
\phi(a)\in\phi(fh)=\phi(f)\phi(h)=1\!\cdot\!\phi(h)=\{\phi(h)\}.
\]
It follows that $\phi(a)=\phi(h)$, so that, by {\color{rot} Lemma 6.2}, $a\in hF$. We, thus, have seen that ${\rm ker}(\phi)$ is a normal closed subset of $H$. \na

\ns\ne Since $1=\phi(1)\in {\rm im}(\phi)$, ${\rm im}(\phi)$ is not empty. \na

\ne Now let $c$ and $d$ be elements in ${\rm im}(\phi)$. Then $H$ contains elements $a$ and $b$ such that $\phi(a)=c$ and $\phi(b)=d$. Thus, by {\color{rot} Lemma 6.1},
\[
c^{*}d=\phi(a)^{*}\phi(b)=\phi(a^{*})\phi(b)=\phi(a^{*}b)\subseteq {\rm im}(\phi),
\]
and, since $c$ and $d$ have been chosen arbitrarily from ${\rm im}(\phi)$, we have shown that ${\rm im}(\phi)$ is a closed subset of $H'$. It remains to show that $H/\!/{\rm ker}(\phi)\cong {\rm im}(\phi)$. \na

\ne We set $F:={\rm ker(\phi)}$. From {\color{rot} (i)} we know that $F$ is normal in $H$. Thus, by {\color{rot} Lemma 6.2}, $\phi(a)=\phi(b)$ if and only if $a^{F}=b^{F}$ for any two elements $a$ and $b$ in $H$. \na

\ne For each element $h$ in $H$, define $\psi(h^{F}):=\phi(h)$. Since $\phi(a)=\phi(b)$ if and only if $a^{F}=b^{F}$ for any two elements $a$ and $b$ in $H$, $\psi$ is an injective map from $H/\!/F$ to ${\rm im}(\phi)$. The definition of $\psi$ also implies that ${\rm im}(\psi)={\rm im}(\phi)$. \na

\ns\ne In order to show that $\psi$ is a homomorphism, we choose elements $a$ and $b$ in $H$. Then
\[
\psi(a^{F}b^{F})=\psi(\{h^{F}\st h\in aFb\})=\{\psi(h^{F})\st h\in aFb\}=\{\phi(h)\st h\in aFb\}=\phi(aFb).
\]
On the other hand, we have
\[
\phi(afb)=\phi(a)\phi(f)\phi(b)=\phi(a)\!\cdot\! 1\!\cdot\!\phi(b)=\phi(a)\phi(b)=\phi(ab)
\]
for each element $f$ in $F$. Thus,
\[
\psi(a^{F}b^{F})=\phi(ab)=\phi(a)\phi(b)=\psi(a^{F})\psi(b^{F}),
\]
and we are done. \qed\gra

\ns\ne {\bf Theorem 6.4} \na

\ne {\it Let $H$ be a hypergroup, let $D$ and $E$ be closed subsets of $H$, and assume that $D\subseteq E$. Assume that $E$ is normal in $H$. Then $E/\!/D$ is a normal closed subset of $H/\!/D$ and $(H/\!/D)/\!/(E/\!/D)\cong H/\!/E$.} \gra

\ne {\it Proof.} That $E/\!/D$ is a normal closed subset of $H/\!/D$ follows from {\color{rot} Lemma 5.3(i)}. It remains to prove that $H/\!/D$ and $(H/\!/D)/\!/(E/\!/D)\cong H/\!/E$. \na

\ne Let $a$ and $b$ be elements in $H$ with $a^{D}=b^{D}$. Then $DaD=DbD$. Thus, as $D\subseteq E$, $EaE=EbE$. It follows that $a^{E}=b^{E}$. This shows that
\[
\phi\!:\ H/\!/D\ \rightarrow\ H/\!/E,\ \ h^{D}\ \mapsto\ h^{E}
\]
is a surjective map. \na

\ne Note that
\[
\phi(a^{D}b^{D})=\phi(\{h^{D}\st h\in aDb\})=\{h^{E}\st h\in aDb\}=aDb/\!/E
\]
and that
\[
\phi(a^{D})\phi(b^{D})=a^{E}b^{E}=\{h^{E}\st h\in aEb\}=aEb/\!/E
\]
\ns\ne for any two elements $a$ and $b$ in $H$. On the other hand, since $E$ is assumed to be normal in $H$, we obtain from {\color{rot} Lemma 4.1(i)} that
\[
aDb/\!/E=EaDbE/\!/E=aEb/\!/E
\]
for any two elements $a$ and $b$ in $H$. Thus, we have
\[
\phi(a^{D}b^{D})=\phi(a^{D})\phi(b^{D})
\]
for any two elements $a$ and $b$ in $H$. Since we also have $\phi(1^{D})=1^{E}$, $\phi$ is a homomorphism from $H/\!/D$ to $H/\!/E$. \na

\ne Note finally that, for each element $h$ in $H$, $\phi(h^{D})=1^{E}$ if and only if $h\in E$. Thus, ${\rm ker}(\phi)=E/\!/D$. It follows that
\[
(H/\!/D)/\!/(E/\!/D)=(H/\!/D)/\!/{\rm ker}(\phi)\cong {\rm im}(\phi)=H/\!/E;
\]
cf.\ {\color{rot} Theorem 6.3}. \qed\gra

\ns\ne {\bf Theorem 6.5} \na

\ne {\it Let $H$ be a hypergroup, let $D$ and $E$ be closed subsets of $H$, and assume that $D$ normalizes $E$. Then $ED$ is a closed subset of $H$, $E$ is a normal closed subset of $ED$, $E\cap D$ is a normal closed subset of $D$, and $ED/\!/E\cong D/\!/(E\cap D)$.} \gra

\ne {\it Proof.} That $ED$ is a closed subset of $H$ was shown in {\color{rot} Lemma 4.1(ii)}, that $E$ is a normal closed subset of $ED$ was shown in {\color{rot} Lemma 4.1(iii)}, and that $E\cap D$ is a normal closed subset of $D$ was shown in {\color{rot} Lemma 4.1(iv)}. It remains to show that $ED/\!/E\cong D/\!/(E\cap D)$. \na

\ne For each element $d$ in $D$, define $\psi(d):=d^{E}$. Then
\[
\psi(ab)=\{\psi(d)\st d\in ab\}=\{d^{E}\st d\in ab\}=ab/\!/E
\]
and
\[
\psi(a)\psi(b)=a^{E}b^{E}=\{d^{E}\st d\in aEb\}=aEb/\!/E
\]
for any two elements $a$ and $b$ in $D$. \na

\ns\ne On the other hand, since we are assuming that $D$ normalizes $E$, we obtain from {\color{rot} Lemma 4.1(i)} that
\[
ab/\!/E=EabE/\!/E=aEb/\!/E
\]
for any two elements $a$ and $b$ in $D$. Thus, we have $\psi(ab)=\psi(a)\psi(b)$ for any two elements $a$ and $b$ in $D$. Since we also have $\psi(1)=1^{E}$, $\psi$ is a homomorphism from $D$ to $ED/\!/E$. \na

\ne Note finally that
\[
{\rm ker}(\psi)=\{d\in D\st d\in E\}=E\cap D.
\]
Thus, by {\color{rot} Theorem 6.3},
\[
D/\!/(E\cap D)=D/\!/{\rm ker}(\psi)\cong {\rm im}(\psi)=ED/\!/E,
\]
and we are done. \qed\gra\gra

\ns\centerline{\bf 7. Solvable hypergroups} \gra\na

\ne A hypergroup $H$ is called {\color{blau} \em solvable} if it contains closed subsets $F_{0}$, $\ldots$, $F_{n}$ such that $F_{0}=\{1\}$, $F_{n}=H$, and, for each element $i$ in $\{0,\ldots,n\}$ with $1\le i$, $F_{i-1}\subseteq F_{i}$, $F_{i}/\!/F_{i-1}$ is thin, and $|F_{i}/\!/F_{i-1}|$ is a prime number. \gra

\ns\ne {\bf Lemma 7.1} \na

\ne {\it Closed subsets of solvable hypergroups are solvable.} \gra

\ne {\it Proof.} Let $H$ be a solvable hypergroup. Then $H$ contains closed subsets $F_{0}$, $\ldots$, $F_{n}$ such that $F_{0}=\{1\}$, $F_{n}=H$, and, for each element $i$ in $\{1,\ldots,n\}$, $F_{i-1}\subseteq F_{i}$, $F_{i}/\!/F_{i-1}$ is thin, and $|F_{i}/\!/F_{i-1}|$ is a prime number. \na

\ne Let $E$ be a closed subset of $H$. For each element $i$ in $\{0,\ldots,n\}$, we set $E_{i}:=F_{i}\cap E$. \na

\ne Let $i$ be an element in $\{1,\ldots,n\}$. Then $E_{i}\subseteq F_{i}$. Thus, as $F_{i}$ normalizes $F_{i-1}$, $E_{i}$ normalizes $F_{i-1}$. It follows that $F_{i-1}E_{i}$ is a closed subset of $F_{i}$; cf.\ {\color{rot} Lemma 4.1(ii)}. Thus, by {\color{rot} Theorem 5.2}, $F_{i-1}E_{i}/\!/F_{i-1}$ is a closed subset of $F_{i}/\!/F_{i-1}$. Since $F_{i}/\!/F_{i-1}$ is thin and $|F_{i}/\!/F_{i-1}|$ is a prime number, this implies that
\[
F_{i-1}/\!/F_{i-1}=F_{i-1}E_{i}/\!/F_{i-1}\qquad {\rm or}\qquad F_{i-1}E_{i}/\!/F_{i-1}=F_{i}/\!/F_{i-1}.
\]
\na

\ns\ne In the first case, we obtain that $F_{i-1}=F_{i-1}E_{i}$, and that implies that $E_{i}\subseteq F_{i-1}$. Now, as $E_{i-1}=F_{i-1}\cap E$, we conclude that $E_{i-1}=E_{i}$. \na

\ne In the second case, we recall that, since $E_{i}$ normalizes $F_{i-1}$,
\[
F_{i-1}E_{i}/\!/F_{i-1}\cong E_{i}/\!/E_{i-1};
\]
cf.\ {\color{rot} Theorem 6.5}. Thus, we obtain in this case that $E_{i}/\!/E_{i-1}\cong F_{i}/\!/F_{i-1}$. Now, as $F_{i}/\!/F_{i-1}$ is thin and $|F_{i}/\!/F_{i-1}|$ is a prime number, $E_{i}/\!/E_{i-1}$ is thin, and $|E_{i}/\!/E_{i-1}|$ is a prime number. Since $E_{0}=\{1\}$ and $E_{n}=E$, this shows that $E$ is solvable. \qed\gra

\ns\ne {\bf Lemma 7.2} \na

\ne {\it Let $H$ be a solvable hypergroup, and let $E$ be a normal closed subset of $H$. Then $H/\!/E$ is solvable.}\gra

\ne {\it Proof.} We are assuming that $H$ is solvable. Thus, $H$ contains closed subsets $F_{0}$, $\ldots$, $F_{n}$ such that $F_{0}=\{1\}$, $F_{n}=H$, and, for each element $i$ in $\{1,\ldots,n\}$, $F_{i-1}\subseteq F_{i}$, $F_{i}/\!/F_{i-1}$ is thin, and $|F_{i}/\!/F_{i-1}|$ is a prime number. \na

\ne Let $i$ be an element in $\{1,\ldots,n\}$. Then $F_{i-1}$ normalizes $E$. Thus, by {\color{rot} Lemma 4.1(ii)}, $EF_{i-1}$ is a closed subset of $H$. It follows that  $EF_{i-1}\cap F_{i}$ is a closed subset of $F_{i}$ and $F_{i-1}\subseteq EF_{i-1}\cap F_{i}$. Since $F_{i}/\!/F_{i-1}$ is thin and $|F_{i}/\!/F_{i-1}|$ is a prime number, this implies that
\[
F_{i-1}=EF_{i-1}\cap F_{i}\qquad {\rm or}\qquad EF_{i-1}\cap F_{i}=F_{i}
\]
\na

\ns\ne On the other hand, since $F_{i}$ normalizes $E$ and $F_{i-1}$, we have $EF_{i-1}h\subseteq EhF_{i-1}\subseteq hEF_{i-1}$ for each element $h$ in $F_{i}$, so $F_{i}$ normalizes $EF_{i-1}$. Thus, by {\color{rot} Theorem 6.5},
\[
EF_{i}/\!/EF_{i-1}\cong F_{i}/\!/(EF_{i-1}\cap F_{i}).
\]
\na

\ne It follows that $EF_{i}/\!/EF_{i-1}=F_{i}/\!/F_{i-1}$ or $EF_{i}/\!/EF_{i-1}=F_{i}/\!/F_{i}$. \na

\ne In the first case, we obtain that $EF_{i}/\!/EF_{i-1}$ is thin and $|EF_{i}/\!/EF_{i-1}|$ is a prime number. From the fact that $EF_{i}/\!/EF_{i-1}$ is thin we obtain that $EF_{i-1}$ is strongly normal in $EF_{i}$; cf.\ {\color{rot} Lemma 5.4}. In particular, $EF_{i-1}$ is normal in $EF_{i}$. Thus, by {\color{rot} Theorem 6.4},
\[
(EF_{i}/\!/E)/\!/(EF_{i-1}/\!/E)\cong EF_{i}/\!/EF_{i-1}.
\]
It follows that $(EF_{i}/\!/E)/\!/(EF_{i-1}/\!/E)$ is thin and $|(EF_{i}/\!/E)/\!/(EF_{i-1}/\!/E)|$ is a prime number. \na

\ne In the second case, $EF_{i-1}=EF_{i}$. \na

\ne Since $EF_{n}/\!/E=H/\!/E$ and $EF_{0}/\!/E=E/\!/E$, this shows that $H/\!/E$ is solvable. \qed\gra

\ns\ne {\bf Lemma 7.3} \na

\ne {\it Let $H$ be a hypergroup, and let $E$ be a closed subset of $H$. Assume that $E$ and $H/\!/E$ are solvable. Then $H$ is solvable.} \gra

\ne {\it Proof.} We are assuming that $E$ is solvable. Thus, $E$ contains closed subsets $F_{0}$, $\ldots$, $F_{m}$ such that $F_{0}=\{1\}$, $F_{m}=E$, and, for each element $i$ in $\{1,\ldots,m\}$, $F_{i-1}\subseteq F_{i}$, $F_{i}/\!/F_{i-1}$ is thin, and $|F_{i}/\!/F_{i-1}|$ is a prime number. \na

\ne We are assuming that $H/\!/E$ is solvable. Thus, by {\color{rot} Theorem 5.2}, $H$ contains closed subsets $F_{m+1}$, $\ldots$, $F_{n}$ such that $F_{n}/\!/E=H/\!/E$ and, for each element $i$ in $\{m+1,\ldots,n\}$, $F_{i-1}\subseteq F_{i}$, $(F_{i}/\!/E)/\!/(F_{i-1}/\!/E)$ is thin, and $|(F_{i}/\!/E)/\!/(F_{i-1}/\!/E)|$ is a prime number. \na

\ns\ne Let $i$ be an element in $\{m+1,\ldots,n\}$. Then $(F_{i}/\!/E)/\!/(F_{i-1}/\!/E)$ is thin. Thus, by {\color{rot} Lemma 5.4}, $F_{i-1}/\!/E$ is strongly normal in $F_{i}/\!/E$. Thus, by {\color{rot} Lemma 5.3(ii)}, $F_{i-1}$ is a strongly normal in $F_{i}$. In particular, $F_{i-1}$ is a normal in $F_{i}$. Thus, by {\color{rot} Theorem 6.4},
\[
(F_{i}/\!/E)/\!/(F_{i-1}/\!/E)\cong F_{i}/\!/F_{i-1}.
\]
Now, as $(F_{i}/\!/E)/\!/(F_{i-1}/\!/E)$ is thin and $|(F_{i}/\!/E)/\!/(F_{i-1}/\!/E)|$ is a prime number, $F_{i}/\!/F_{i-1}$ is thin and $|F_{i}/\!/F_{i-1}|$ is a prime number. \na

\ne What we have seen so far is that, for each element $i$ in $\{1,\ldots,n\}$, $F_{i-1}\subseteq F_{i}$, $F_{i}/\!/F_{i-1}$ is thin, and $|F_{i}/\!/F_{i-1}|$ is a prime number. From $F_{n}/\!/E=H/\!/E$ we also obtain that $F_{n}=H$. Thus, as $F_{0}=\{1\}$, $H$ is solvable. \qed\gra

\ns\ne Let $H$ be a hypergroup. Recall from {\color{rot} Section 4} that a closed subset $E$ of $H$ is said to be subnormal in $H$ if $H$ contains closed subsets $F_{0}$, $\ldots$, $F_{n}$ such that $F_{0}=E$, $F_{n}=H$, and, for each element $i$ in $\{1,\ldots,n\}$, $F_{i-1}$ is normal in $F_{i}$. \na

\ns\ne We are now in the position to weaken the hypothesis of {\color{rot} Lemma 7.2}. \gra

\ns\ne {\bf Theorem 7.4} \na

\ne {\it Let $H$ be a solvable hypergroup, and let $E$ be a subnormal closed subset of $H$. Then $H/\!/E$ is solvable.} \gra

\ne {\it Proof.} Clearly, if $E=H$, $H/\!/E$ is solvable. Therefore, we assume that $E\neq H$. \na

\ne From $E\neq H$ together with the fact that $E$ is a subnormal closed subset of $H$ we obtain a proper normal closed subset $F$ of $H$ such that $E$ is subnormal in $F$. By induction, $F/\!/E$ is solvable. \na

\ne Since $H$ is solvable and $F$ is a normal closed subset of $H$, $H/\!/F$ is solvable; cf.\ {\color{rot} Lemma 7.2}. On the other hand, as $F$ is a normal closed subset of $H$, we obtain from {\color{rot} Theorem 6.4} also that
\[
(H/\!/E)/\!/(F/\!/E)\cong H/\!/F.
\]
Thus, $(H/\!/E)/\!/(F/\!/E)$ is solvable. \na

\ne Now, as both $F/\!/E$ and $(H/\!/E)/\!/(F/\!/E)$ are solvable, so is $H/\!/E$; cf. {\color{rot} Lemma 7.3}. \qed\gra

\ns\ne {\bf Corollary 7.5} \na

\ne {\it Let $H$ be a solvable hypergroup, and let $D$ and $E$ be closed subsets of $H$ with $D\subseteq E$. Assume that $D$ is subnormal in $H$, and that $E/\!/D$ is subnormal in $H/\!/D$. Then $E$ is subnormal in $H$.} \gra

\ne {\it Proof.} We are assuming that $H$ is solvable and that $D$ is subnormal in $H$. Thus, by {\color{rot} Theorem 7.4}, $H/\!/D$ is solvable. Since $H/\!/D$ is solvable and $E/\!/D$ is assumed to be subnormal in $H/\!/D$, $(H/\!/D)/\!/(E/\!/D)$ is solvable; cf.\ {\color{rot} Theorem 7.4}. Thus, by {\color{rot} Theorem 5.2}, $H$ contains closed subsets $F_{0}$, $\ldots$, $F_{n}$ such that $F_{0}/\!/D=E/\!/D$, $F_{n}/\!/D=H/\!/D$, and, for each element $i$ in $\{1,\ldots,n\}$, $F_{i-1}/\!/D\subseteq F_{i}/\!/D$ and $(F_{i}/\!/D)/\!/(F_{i-1}/\!/D)$ is thin. \na

\ns\ne Let $i$ be an element in $\{1,\ldots,n\}$. Since $(F_{i}/\!/D)/\!/(F_{i-1}/\!/D)$ is thin, $F_{i-1}/\!/D$ is strongly normal in $F_{i}/\!/D$; cf.\ {\color{rot} Lemma 5.4}. Thus, by {\color{rot} Lemma 5.3(ii)}, $F_{i-1}$ is strongly normal in $F_{i}$. In particular, $F_{i-1}$ is normal in $F_{i}$; cf.\ {\color{rot} Lemma 4.3}. \na

\ne From $E/\!/D=F_{0}/\!/D$ we obtain that $E=F_{0}$, from $F_{n}/\!/D=H/\!/D$ that $F_{n}=H$. Also, for each element $i$ in $\{1,\ldots,n\}$, $F_{i-1}$ is normal in $F_{i}$. Thus, $E$ is subnormal in $H$. \qed\gra\gra

\ns\centerline{\bf 8. Hall subsets of solvable association schemes} \gra\na

\ne In this section, we turn to association schemes. We will refer freely to some notation introduced in {\color{gruen} [8]}. \na

\ne Our first lemma is a slight generalization of {\color{gruen} [4; Lemma 5.2]}. \gra

\ns\ne {\bf Lemma 8.1} \na

\ne {\it Let $S$ be an association scheme, let $T$ be a closed subset of $S$, let $s$ be an element of $S$, and let $\pi$ be a set of prime numbers. Assume that $s$ is  $\pi$-valenced and that $n_{T}$ is a $\pi$-number. Then $s^{T}$ is $\pi$-valenced.} \gra

\ne {\it Proof.} By {\color{gruen} [8; Theorem 4.1.3(ii)]}, $n_{s^{T}}n_{T}=n_{TsT}$. Moreover, by {\color{gruen} [8; Lemma 2.3.3]}, $n_{TsT}$ divides $n_{T}n_{s}n_{T}$. Thus, $n_{s^{T}}$ divides $n_{T}n_{s}$. Thus, as $n_{T}$ and $n_{s}$ are assumed to be $\pi$-numbers, $n_{s^{T}}$ is a $\pi$-number, and that means that $s^{T}$ is $\pi$-valenced. \qed\gra

\ns\ne {\bf Lemma 8.2} \na

\ne {\it Let $\pi$ be a set of prime numbers, and let $S$ be a solvable and $\pi$-valenced association scheme. Assume that $S$ has no thin subnormal closed subset different from $\{1\}$ the valency of which is a $\pi$-number. Then $S$ is thin.} \gra

\ne {\it Proof.} Since $S$ is assumed to be solvable, $S$ contains closed subsets $T_{0}$, $\ldots$, $T_{n}$ such that $T_{0}=\{1\}$, $T_{n}=S$, and, for each element $i$ in $\{1,\ldots,n\}$, $T_{i-1}$ is strongly normal in $T_{i}$. \na 

\ne Assume, by way of contradiction, that $S$ is not thin. Then $\{1,\ldots,n\}$ contains an element $i$ such that $T_{i-1}$ is thin and $T_{i}$ is not thin. Since $T_{i-1}$ is strongly normal in $T_{i}$, $O^{\vartheta}(T_{i})\subseteq T_{i-1}$. Thus, as $T_{i-1}$ is thin, $O^{\vartheta}(T_{i})$ is thin. \na

\ns\ne Let $s$ be an element in $T_{i}$. Then, as $O^{\vartheta}(T_{i})$ is thin, $s^{*}s$ is a thin closed subset of $S$; cf.\ {\color{gruen} [8; Lemma 6.7.1(i), (iv)]}. Furthermore, by {\color{gruen} [8; Lemma 6.7.1(v)]}, $s^{*}s$ is strongly normal in $O^{\vartheta}(T_{i})$, and, by {\color{gruen} [8; Lemma 6.7.1(iii)]}, $\{s\}=ss^{*}s$. This latter equation yields $n_{s^{*}s}=n_{s}$; cf.\  {\color{gruen} [8; Lemma 1.4.4(ii)]}. Thus, as $S$ is assumed to be $\pi$-valenced, $n_{s^{*}s}$ is a $\pi$-number. \na

\ne What we have seen is that $s^{*}s$ is a thin subnormal closed subset of $S$ the valency of which is a $\pi$-number. Since $S$ is assumed to have no thin subnormal closed subset different from $\{1\}$ the valency of which is a $\pi$-number, $s^{*}s=\{1\}$. This implies that $s$ is thin. \na

\ne Since $s$ has been chosen arbitrarily from $T_{i}$, we have shown that $T_{i}$ is thin, contradiction. \qed\gra

\ns\ne Let $S$ be an association scheme, and recall from {\color{rot} Section 1} that $S$ is a hypergroup with respect to the hypermultiplication defined by the complex multiplication in $S$. Note that the closed subsets of $S$ are exactly the closed subsets of this hypergroup. Furthermore, the strongly normal closed subsets of $S$ are exactly the strongly normal closed subsets of this hypergroup, so that, by {\color{gruen} [9; Lemma 4.2.5(ii)]}, $S$ is solvable if and only if this hypergroup is solvable. We shall now see that the transition from scheme theory to hypergroup theory is also compatible with quotients. \gra

\ns\ne {\bf Lemma 8.3} \na

\ne {\it Let $S$ be an association scheme, and let $T$ be a closed subset of $S$. Then the hypergroup defined by the complex multiplication in the quotient scheme $S/\!/T$ is equal to the quotient of the hypergroup defined by the complex multiplication in the scheme $S$ over $T$.} \gra

\ne {\it Proof.} By {\color{gruen} [8; Lemma 4.1.4]}, the map which associates to each element $s^{T}$ in $S/\!/T$ the double coset $TsT$ is a bijective map from the quotient scheme $S/\!/T$ to the set of all double cosets of $T$ in $S$. The elements of the quotient of the hypergroup defined by the complex multiplication in the scheme $S$ over $T$ are the double cosets of $T$ in the hypergroup $S$. \na

\ne Let $p$ and $q$ be elements in $S$. Then, by {\color{gruen} [8; Lemma 4.1.4]},
\[
p^{T}q^{T}:=\{s^{T}\st s\in TpTqT\}.
\]
The product $p^{T}$ and $q^{T}$ in the quotient of the hypergroup defined by the complex multiplication in the scheme $S$ over $T$ is defined to be the set $\{s^{T}\st s\in pTq\}$. On the other hand, by {\color{rot} Lemma 5.1}, $\{s^{T}\st s\in pTq\}=\{s^{T}\st s\in TpTqT\}$. \qed\gra

\ns\ne A closed subset $T$ of an association scheme $S$ is said to be {\color{blau} \em subnormal in} $S$ if $T$ is subnormal in the hypergroup defined by the complex multiplication in $S$. \gra

\ns\ne {\bf Lemma 8.4} \na

\ne {\it Let $S$ be a solvable association scheme, and let $T$ be a subnormal closed subset of $S$. Then we have the following.} \tabulaturna

\ne\loss

\hfill (i) & {\it The quotient scheme $S/\!/T$ is solvable.} \\

\end{tabular} \tabulaturkla

\ne\loss

\hfill (ii) & {\it Let $U$ be a closed subset of $S$ with $T\subseteq U$. Assume that $U/\!/T$ is subnormal in $S/\!/T$. Then $U$ is subnormal in $S$.} \\

\end{tabular} \gra

\ne {\it Proof.} (i) We are assuming that $T$ is a subnormal closed subset of $S$. Thus, $T$ is a subnormal closed subset of the hypergroup $S$. Thus, by {\color{rot} Theorem 7.4}, the quotient $S/\!/T$ is solvable. Thus, by {\color{rot} Lemma 8.3}, the quotient scheme $S/\!/T$ is solvable. \na

\ne (ii) We are assuming that $T$ is a subnormal closed subset of the scheme $S$. Thus, $T$ is a subnormal closed subset of the hypergroup $S$ defined by the complex multiplication in $S$. \na

\ne We are assuming that the quotient scheme $U/\!/T$ is subnormal in the quotient scheme $S/\!/T$. Thus, by {\color{rot} Lemma 8.3}, the quotient $U/\!/T$ is subnormal in the quotient $S/\!/T$. Thus, by {\color{rot} Corollary 7.5}, $U$ is subnormal in the hypergroup $S$ defined by its complex multiplication. Thus, $U$ is subnormal in the scheme $S$. \qed\gra

\ns\ne {\bf Theorem 8.5} \na

\ne {\it Let $\pi$ be a set of prime numbers, and let $S$ be a solvable and $\pi$-valenced association scheme. Then $S$ contains a strongly normal closed $\pi$-subset which contains all subnormal closed $\pi$-subsets of $S$.} \gra

\ne {\it Proof.} Let $U$ be maximal among the subnormal closed $\pi$-subsets of $S$. By hypothesis, $S$ is solvable. Thus, as $U$ is a subnormal closed subset of $S$, $S/\!/U$ is solvable; cf.\ {\color{rot} Lemma 8.4(i)}. Furthermore, since $S$ is $\pi$-valenced and $n_{U}$ is a $\pi$-number, we obtain from {\color{rot} Lemma 8.1} that $S/\!/U$ is $\pi$-valenced. Thus, we may apply {\color{rot} Lemma 8.2} to the quotient scheme $S/\!/U$ in place of $S$. \na

\ne Let $V$ be a closed subset of $S$ with $U\subseteq V$. Assume that $V/\!/U$ is thin, that $V/\!/U$ is subnormal in $S/\!/U$, and that $n_{V/\!/U}$ is a $\pi$-number. \na

\ns\ne Since $U$ is subnormal in $S$ and $V/\!/U$ is subnormal in $S/\!/U$, we obtain from {\color{rot} Lemma 8.4(ii)} that $V$ is subnormal in $S$. On the other hand, $n_{U}$ as well as $n_{V/\!/U}$ are $\pi$-numbers. Thus, by {\color{gruen} [8; Lemma 4.3.3(i)]}, $n_{V}$ is a $\pi$-number, so that the choice of $U$ forces $U=V$. \na

\ne This shows that $S/\!/U$ does not contain a thin subnormal closed subset different from $U/\!/U$ the valency of which is a $\pi$-number. Thus, by {\color{rot} Lemma 8.2}, $S/\!/U$ is thin. It follows that $U$ is strongly normal in $S$; cf.\ {\color{gruen} [8; Lemma 4.2.5(ii)]}. \na

\ne Let $T$ be a subnormal closed $\pi$-subset of $S$. Since $U$ is strongly normal in $S$, $U$ is normal in $S$; cf.\ {\color{gruen} [8; Lemma 2.5.5]}. Thus, by {\color{gruen} [8; Lemma 2.5.2(iii)]}, $UT$ is a closed subset of $S$ and is subnormal in $S$. \na

\ne On the other hand, by {\color{gruen} [8; Lemma 2.3.6(i)]}, $n_{U}n_{T}=n_{UT}n_{U\cap T}$. Thus, $UT$ is a closed $\pi$-subset of $S$. Thus, as $U\subseteq UT$, the choice of $U$ forces $T\subseteq U$. \qed\gra

\ns\ne Let $\pi$ be a set of prime numbers, and let $S$ be a solvable and $\pi$-valenced association scheme. From {\color{rot} Theorem 8.5} we obtain that $S$ contains a strongly normal closed $\pi$-subset which contains all subnormal closed $\pi$-subsets of $S$. In the following, this closed subset of $S$ will be denoted by $O_{\pi}(S)$. \gra

\ns\ne {\bf Proposition 8.6} \na

\ne {\it Let $\pi$ be a set of prime numbers, and let $S$ be a solvable and $\pi$-valenced association scheme. Set $O:=O_{\pi}(S)$. Then we have the following.} \tabulaturna

\ne\loss

\hfill (i) & {\it The scheme $S/\!/O$ possesses at least one Hall $\pi$-subset.} \\

\end{tabular} \tabulaturkla

\ne\loss

\hfill (ii) & {\it Any two Hall $\pi$-subsets of $S/\!/O$ are conjugate in $S$.} \\

\end{tabular} \tabulaturkla

\ne\loss

\hfill (iii) & {\it Any closed $\pi$-subset of $S/\!/O$ is contained in a Hall $\pi$-subset of $S/\!/O$.} \\

\end{tabular} \gra

\ne {\it Proof.} Since $O$ is a strongly normal closed subset of $S$, $S/\!/O$ is thin; cf.\ {\color{gruen} [8; Lemma 4.2.5(ii)]}. Moreover, by {\color{rot} Lemma 8.4(i)}, $S/\!/O$ is solvable. Thus, with a reference to the group correspondence, the claim follows from {\color{gruen} [2; Theorem]}. \qed\gra

\ns\ne Now we are in the position to prove the main result of this article. We split the proof into three parts. \gra

\ns\ne {\bf Theorem 8.7} \na

\ne {\it Let $\pi$ be a set of prime numbers, and let $S$ be a solvable and $\pi$-valenced association scheme. Then $S$ possesses at least one Hall $\pi$-subset.} \gra

\ne {\it Proof.} Set $O:=O_{\pi}(S)$. By {\color{rot} Proposition 8.6(i)}, the scheme $S/\!/O$ possesses at least one Hall $\pi$-subset. Let $T$ be a closed subset of $S$ such that $T/\!/O$ is a Hall $\pi$-subset of $S/\!/O$; cf.\ {\color{gruen} [8; Lemma 4.1.7(ii)]}. We will see that $T$ is a Hall $\pi$-subset of $S$. \na

\ne Since $T/\!/O$ is a Hall $\pi$-subset of $S/\!/O$, $n_{T/\!/O}$ is a $\pi$-number. Since $n_{O}$ is a $\pi$-number, too, we obtain from {\color{gruen} [8; Lemma 4.3.3(i)]} that $n_{T}$ is a $\pi$-number. But $T$ is also $\pi$-valenced, since $S$ is $\pi$-valenced. Thus, $T$ is a closed $\pi$-subset of $S$, and it remains to be shown that $n_{S/\!/T}$ is a $\pi'$-number. \na

\ne Since $T/\!/O$ is a Hall $\pi$-subset of $S/\!/O$, $n_{(S/\!/O)/\!/(T/\!/O)}$ is a $\pi'$-number. From {\color{gruen} [8; Lemma 4.3.3]} we also obtain that $n_{S/\!/T}=n_{(S/\!/O)/\!/(T/\!/O)}$. Thus, $n_{S/\!/T}$ is a $\pi'$-number. \qed\gra

\ns\ne {\bf Theorem 8.8} \na

\ne {\it Let $\pi$ be a set of prime numbers, and let $S$ be a solvable and $\pi$-valenced association scheme. Then any two Hall $\pi$-subsets of $S$ are conjugate in $S$.} \gra

\ne {\it Proof.} Let $T$ and $U$ be Hall $\pi$-subsets of $S$, and set $O:=O_{\pi}(S)$. Since $S$ is assumed to be $\pi$-valenced, so is $T$. Thus, as $n_{O}$ is a $\pi$-number, $T/\!/O$ is $\pi$-valenced; cf.\ {\color{rot} Lemma 8.1}. \na

\ne Since $n_{T}$ is a $\pi$-number, so is $n_{T/\!/O}$; cf.\ {\color{gruen} [8; Lemma 4.3.3(i)]}. Thus, $T/\!/O$ is a closed $\pi$-subset of $S/\!/O$. From {\color{gruen} [8; Lemma 4.3.3]} we also obtain that $n_{S/\!/T}=n_{(S/\!/O)/\!/(T/\!/O)}$. Thus, as $n_{S/\!/T}$ is a $\pi'$-number, so is $n_{(S/\!/O)/\!/(T/\!/O)}$. It follows that $T/\!/O$ is a Hall $\pi$-subset of $S/\!/O$. \na

\ne Similarly, $U/\!/O$ is a Hall $\pi$-subset of $S/\!/O$. Thus, by {\color{rot} Proposition 8.6(ii)}, $S$ contains an element $s$ such that $(s^{O})^{*}(T/\!/O)s^{O}=U/\!/O$. Thus, by {\color{gruen} [8; Lemma 4.1.4]}, $s^{*}Ts=U$. \qed\gra

\ns\ne {\bf Theorem 8.9} \na

\ne {\it Let $\pi$ be a set of prime numbers, and let $S$ be a solvable and $\pi$-valenced association scheme. Then any closed $\pi$-subset of $S$ is contained in a Hall $\pi$-subset of $S$.} \gra

\ne {\it Proof.} Let $T$ be a closed $\pi$-subset of $S$, and set $O:=O_{\pi}(S)$. Then, $OT$ is a closed subset of $S$. Thus, by {\color{gruen} [8; Lemma 4.1.7(i)]}, $OT/\!/O$ is a closed subset of $S/\!/O$. \na

\ne Since $S$ is assumed to be $\pi$-valenced, so is $OT$. Since $T$ is a closed $\pi$-subset of $S$, $n_{T}$ is a $\pi$-number. Also $n_{O}$ is a $\pi$-number. Thus, by {\color{gruen} [8; Lemma 4.3.4]}, $n_{OT/\!/O}$ is a $\pi$-number. Thus, $OT/\!/O$ is a closed $\pi$-subset of $S/\!/O$. Thus, by {\color{rot} Proposition 8.6(iii)}, $S/\!/O$ contains a Hall $\pi$-subset which contains $OT/\!/O$. Let $U$ be a closed subset of $S$ such that $O\subseteq U$, $U/\!/O$ is a Hall $\pi$-subset of $S/\!/O$, and $OT/\!/O\subseteq U/\!/O$; cf.\ {\color{gruen} [8; Lemma 4.1.7(ii)]}. \na

\ne We will see that $U$ is a Hall $\pi$-subset of $S$ which contains $T$. \na

\ne Since $U/\!/O$ is a Hall $\pi$-subset of $S/\!/O$, $n_{U/\!/O}$ is a $\pi$-number. Thus, as $n_{O}$ is a $\pi$-number, $n_{U}$ is a $\pi$-number; cf.\ {\color{gruen} [8; Lemma 4.3.3(i)]}. Since $U/\!/O$ is a Hall $\pi$-subset of $S/\!/O$, $n_{(S/\!/O)/\!/(U/\!/O)}$ is a $\pi'$-number. Furthermore, by {\color{gruen} [8; Lemma 4.3.3]}, $n_{S/\!/U}=n_{(S/\!/O)/\!/(U/\!/O)}$. Thus, $n_{S/\!/U}$ is a $\pi'$-number. It follows that $U$ is a Hall $\pi$-subset of $S$. \na

\ne From $OT/\!/O\subseteq U/\!/O$ we obtain that $OT\subseteq U$; cf.\ {\color{gruen} [8; Lemma 4.1.4]}. Thus, $U$ is a Hall $\pi$-subset of $S$, and $U$ contains $T$. \qed\gra\gra

\ns\centerline{\bf References} \gra\na

\ne 1. French, C., Zieschang, P.-H.: A Schur-Zassenhaus Theorem for association schemes, {\it J.\ Algebra} {\bf 435}, 88-123 (2015) \na

\ne 2. Hall, P.: A note on soluble groups, {\it J. London Math. Soc.}\ {\bf 3}, 98-105 (1928) \na

\ne 3. Hanaki, A., Miyamoto, I.: Classification of association schemes of small order, Online catalogue. http://kissme.shinshu-u.ac.jp/as \na

\ne 4. Hirasaka M., Muzychuk, M., Zieschang, P.-H.: A generalization of Sylow's theorems on finite groups to association schemes, {\it Math.\ Z.}\ {\bf 241}, 665-672 (2002) \na

\ne 5. Marty, F.: Sur une g${\rm\acute{e}}$n${\rm\acute{e}}$ralisation de la notion de groupe, in: Huiti${\rm\grave{e}}$me Congres des Math${\rm\acute{e}}$maticiens, Stockholm 1934, 45-59 \na

\ne 6. Rassy, M., Zieschang, P.-H.: Basic structure theory of association schemes, {\it Math.\ Z.}\ {\bf 227}, 391-402 (1998) \na

\ne 7. Tanaka, R., Zieschang, P.-H.: On a class of wreath products of hypergroups and association schemes, {\it J.\ Algebraic Combin.}\ {\bf 37}, 601-619 (2013) \na

\ne 8. Zieschang, P.-H.: {\it Theory of Association Schemes}, Springer Monographs in Mathematics, Berlin Heidelberg New York, 2005 \na

\ne 9. Zieschang, P.-H.: Hypergroups all non-identity elements of which are involutions, in: {\it Advances in algebra} 305-322, Springer Proc.\ Math.\ Stat.\ {\bf 277}, Springer, Cham, 2019

\end{document}